\documentclass[12pt]{article}

\usepackage[T2A]{fontenc}
\usepackage[cp866]{inputenc}

\usepackage{amsfonts}
\usepackage{amssymb}
\usepackage{amsmath}
\usepackage{amsthm}

\def \le {\leqslant}
\def \ge {\geqslant}

\begin{document}

\centerline{\bf On the derivative of the Minkowski question mark function $?(x)$} \vskip+0.5cm
\centerline{\bf Anna A. Dushistova, Nikolai G.
Moshchevitin \footnote{ Research is supported by grants RFFI 06-01-00518, MD-3003.2006.1, NSh-1312.2006.1 and INTAS 03-51-5070
 } } \vskip+0.5cm
\vskip+1.0cm

\centerline{\bf Abstract}

Let $ x = [0;a_1,a_2,...]$ be  the decomposition of the irrational number $x
\in [0,1]$ into regular continued fraction. Then for the derivative of the
Minkowski function $?(x)$ we prove that $?'(x) = +\infty$ provided $
\limsup_{t\to \infty}\frac{a_1+...+a_t}{t} <\kappa_1 =\frac{2\log
\lambda_1}{\log 2} = 1.388^+$, and $?'(x) = 0$ provided $ \liminf_{t\to
\infty}\frac{a_1+...+a_t}{t} >\kappa_2 =  \frac{4L_5-5L_4}{L_5-L_4}= 4.401^+$
(here $ L_j = \log \left(\frac{j+\sqrt{j^2+4}}{2}\right) - j\cdot\frac{\log
2}{2}$). Constants $\kappa_1,\kappa_2$ are the best possible. Also we prove
that $?'(x) = +\infty$ holds for all $x$ with partial quotients bounded by 4.
 \vskip+1.0cm

{\bf 1.\,\,\, The Minkowski function  $?(x)$.} \,\,\, The function $?(x)$ is
defined as follows.
 $ ?(0)=0, ?(1)= 1$,
 if $?(x)$ is defined for successive  Farey fractions
$\frac{p}{q}, \frac{p'}{q'}$ then
$$?\left(\frac{p+p'}{q+q'}\right) = \frac{1}{2} \left( ?\left(\frac{p}{q}\right)+
?\left(\frac{p'}{q'}\right)\right);
$$
for irrational $x$ function $?(x)$ is defined by continuous arguments. This function firstly was considered by H. Minkowski  (see. \cite{MI},
p.p. 50-51) in 1904. $?(x)$ is a continuous increasing function. It has derivative almost everywhere. It satisfies Lipschitz condition  \cite{SA},
\cite{KI}. It is a well-known fact that the derivative  $?'(x)$ can take only two values -  $0$  or $+\infty$. Almost everywhere we have
 $?'(x)=0$.
Also if irrational $x=[0;a_1,...,a_t,...] $ is represented as a regular  continued fraction with natural partial quotients then
$$
?(x) = \frac{1}{2^{a_1-1}} - \frac{1}{2^{a_1+a_2-1}}+ ...+ \frac{(-1)^{n+1}}{2^{a_1+...+a_n-1}}+ ...  .
$$
These and some other results one can find for example in papers \cite{PARA},\cite{PARA2},\cite{SA}.

Here we should note the connection between function $?(x)$ and Stern-Brocot sequences. We remind the reader the definition of Stern-Brocot
sequences $ F_n $, $ n=0, 1, 2, \dots $. First of all let us put
  $ F_0=\{0, 1\}=\{\frac{0}{1}, \frac{1}{1}\} $.
Then for the sequence $ F_n $ treated as increasing sequence of rationals $
 0=x_{0, n}<x_{1, n}< \dots <x_{N \left( n
  \right), n}=1, N(n)=2^{n},
$  $x_{j, n} =p_{j, n}/q_{j, n}, ( p_{j, n},q_{j, n}) = 1$ we define the next sequence $F_{n+1}$ as  $ F_{n+1} = F_n \cup Q_{n+1} $ where  $
Q_{n+1} $ is the set of the form $ Q_{n+1}=\{ x_{j-1,n}\oplus x_{j,n}, i=1, \dots , N(n)\}.  $ Here operation $\oplus$ means taking the mediant
fraction for  two rational numbers:
  $ \frac{a}{b}\oplus \frac{c}{d} = \frac{a+c}{b+d}$.
The Minkowski question mark function $?(x)$ is the limit distribution function for Stern-Brocot sequences:
$$
?(x)  =\lim_{n\to \infty} \frac{\# \{ \xi\in F_n :\,\,\,\xi \le x\}}{2^n +1}.
$$

 {\bf 2.\,\,\, Notation and parameters.} \,\,\, In this paper $ [0;a_1,...,a_t,...] $
denotes a regular continued fraction with natural partial quotients $a_t$.  $
k_t(a_1,...,a_t) $ denotes continuant. For a continued fraction under
consideration the convergent fraction of order  $t$ is denoted as $ p_t/q_t =
[0;a_1,...,a_t]$ (hence, $ q_t =k_t(a_1,...,a_t) $).
 We need numbers
$$
\lambda_j = \frac{j+\sqrt{j^2+4}}{2}, \,\,\, L_j = \log \lambda_j - j\cdot\frac{\log 2}{2} .$$ Here $j<\lambda_j < j+1$. Note that
\begin{equation}
L_2>L_3>L_1>L_4>0>L_5>L_6> ... \label{123}
\end{equation}
 and
\begin{equation}
  \frac{L_5}{L_5- L_4} \ge \frac{1}{2}.
  \label{12}
  \end{equation}
 Also we need the
  values of continuants
$$
k_{l,j} = k_l(\underbrace{j,..., j}_{l}),\,\,\, k_{0,j} =1,\,\,\, k_{1,j}=j .
$$
From recursion $k_{l+1,j}=jk_{l,j}+k_{l-1,j}$ we deduce
$$
k_{l,j} = c_{1,j} \lambda_j^l+ c_{2,j}(-\lambda_j)^{-l}
$$
where
$$
 c_{1,j} + c_{2,j} =1,\,\,\,
 c_{1,j} \lambda_j - c_{2,j}(\lambda_j)^{-1} =j.
 $$
 Hence
 $$
 1-\frac{j}{j^2+1}
 <c_{1,j} <1,\,\,\,
 0< c_{2,j} <
\frac{j}{j^2+1}
$$
and
\begin{equation}
k_{l,j} <   \lambda_j^l.  \label{kont}
\end{equation}
  Also
we should consider the constants
\begin{equation}
\kappa_1 = \frac{2\log \lambda_1}{\log 2} = 1.388^+,\,\,\,\, \kappa_2 = \frac{4L_5-5L_4}{L_5-L_4}= 4.401^+. \label{kappa}
\end{equation}
For a natural  $n$ and a $n$-tuple of nonnegative integer numbers $(r_1,...,r_n)$ we put $
 t = \sum_{j=1}^nr_j.
$ Now we define the set
$$
W_n(r_1,...,r_n) = \{ (a_1,...,a_t)\,\,:\,\,\, \#\{i\,:\, a_i = j\} = r_j\} .$$
Let
\begin{equation}
\mu_n (r_1,...,r_n) = \max_{ (a_1,...,a_t) \in W_n(r_1,...,r_n)} k_t (a_1,...,a_t). \label{MU}
\end{equation}
For real positive  $\omega$ we define
$$
\Omega_{\omega ,n, t} = \left\{ (r_1,...,r_n)\, :\,\, r_j \in \mathbb{N}_0,\,\, \sum_{j=1}^n (j-\omega )r_j \ge0,\,\,
  \sum_{j=1}^nr_j
=t \right\} .
$$
Let $\omega = \kappa_2 +\eta  < 5$ and $\eta \in [ 0, 1/2)$. It is easy to see
that for any  $ n \ge 5$ the following unequality is valid:
\begin{equation}
\max_{(r_1,...,r_n)\in \Omega_{\kappa_2+\eta, n, t} } \sum_{j=1}^n r_jL_j 
\leq
(L_5-L_4)t\eta ,\,\,\,\, L_5 - L_4 < 0. \label{MAIN}
\end{equation}
We give the proof of  (\ref{MAIN}) in section 5.

Also for  $ r_1 \ge 1$  we consider the set
$$
V_n(r_1,...,r_n) = \{ (a_1,...,a_t)\,\,:\,\,\, \#\{i\,:\, a_i = j\} = r_j, \,\,\, a_1 = 1\}.
$$
Let
$$ k[r_1,...,r_n] = k_t( \underbrace{1,..., 1}_{r_1},  \underbrace{2,..., 2}_{r_2},..., \underbrace{n,..., n}_{r_n}).
$$

 I.D. Kan in \cite{KAN} proved the following statement.

{\bf Lemma 1.}
$$
 \max_{ (a_1,...,a_t) \in V_n(r_1,...,r_n)} k_t (a_1,...,a_t)
= k[r_1,...,r_n].
$$

We should note that Lemma 1 is a generalization of a result from  \cite{MOZ}.

To get an upper bound for $k[r_1,...,r_n]$ we use formula
\begin{equation}
k_{t+l} (a_1,...,a_t, b_1,...,b_l) = k_t (a_1,...,a_t)  k_l (b_1,...,b_l) +
 k_{t-1} (a_1,...,a_{t-1})  k_{l-1} (b_2,...,b_l).
 \label{KNUTH}
 \end{equation}
Let $r_{h_1},...,r_{h_f},\,\,\, 1\le h_1<...< h_f =n$
 be all {\it positive} numbers from the set $r_1,...,r_n$.
 Here $h_j \ge j$.
 Then from (\ref{KNUTH})
 and  inequalities
 $$
k_{r_{h_{j+1}}-1,h_{j+1}}\le k_{r_{h_{j+1}},h_{j+1}}/h_{j+1},\,\,\, k[r_1,...,r_{h_j}-1]\le k[r_1,...,r_{h_j}]/h_j
$$
  we deduce the inequality
$$
k[r_1,...,r_{h_j},\underbrace{0,..., 0}_{h_{j+1}-h_j-1}, r_{h_{j+1}}] =k[r_1,...,r_{h_j}]k_{r_{h_{j+1}},h_{j+1}} + k[r_1,...,r_{h_j}-1]
k_{r_{h_{j+1}}-1,h_{j+1}}\le
$$
$$\le
k[r_1,...,r_{h_j}]k_{r_{h_{j+1}},h_{j+1}} \left(1+\frac{1}{h_jh_{j+1}}\right).
$$
Now
\begin{equation}
k[r_1,...,r_n]\le \prod_{j=1}^n k_{r_j,j} \prod_{j=1}^{f-1}\left(1+\frac{1}{h_jh_{j+1}}\right) \le \prod_{j=1}^n k_{r_j,j}
\prod_{j=1}^{n-1}\left(1+\frac{1}{j(j+1)}\right). \label{prod}
\end{equation}
But
$$
\prod_{j=1}^{n-1}\left(1+\frac{1}{j(j+1)}\right)\le \prod_{j=1}^{+\infty}\left(1+\frac{1}{j(j+1)}\right) \le e.
$$
 Hence from
 Lemma  1, inequalities (\ref{prod},\ref{kont})
and
$$
k_t(a_1,...,a_t) \le k_{t+1}(1,a_1,...,a_t)
$$
as a corollary we deduce the following upper bound for $\mu_n (r)$:
 \begin{equation}
\mu_n (r_1,...,r_n) \le   \lambda_1 e \prod_{j=1}^n \lambda_j^{r_j}. \label{PRT}
\end{equation}

{\bf 3.\,\,\, A result by J. Paradis,  P. Viader, L. Bibiloni.}\,\,\, In
\cite{PARA2} the following statement is proved.

{\bf Theorem A.}\,\,\,{\it

1. Let for real irrational $x\in (0,1)$ in continued fraction expansion $
x=[0;a_1,...,a_t,...] $ with $\kappa_1 $ from (\ref{kappa}) the following
inequality be valid:
$$
 \limsup_{t\to
\infty}\frac{a_1+...+a_t}{t} <\kappa_1 .$$ Then if
 $?'(x)$ exists the equality $?'(x)=+\infty $  holds.

 2. Let $\kappa_3 = 5.319^+$
be the root of equation
 $\frac{2\log (1+x)}{\log 2} - x =0$.
 Let for real irrational $x\in (0,1)$ in continued fraction expansion $
x=[0;a_1,...,a_t,...] $ holds
$$
 \liminf_{t\to
\infty}\frac{a_1+...+a_t}{t} \ge \kappa_3.$$ Then if
 $?'(x)$ exists the equality $?'(x)=0$  holds.
 }

{\bf 4.\,\,\, New results.}\,\,\, Here we give the stronger version of the Theorem
A.

{\bf Theorem 1.}\,\,\,{\it

1. Let for real irrational $x\in (0,1)$ in continued fraction expansion $
x=[0;a_1,...,a_t,...] $ with $\kappa_1 $ from (\ref{kappa}) the following
inequality be valid:
$$
 \limsup_{t\to
\infty}\frac{a_1+...+a_t}{t} <\kappa_1 .$$ Then
  $?'(x)$ exists and $?'(x)=+\infty $.

2. For any positive $\varepsilon$ there exists a quadratic irrationality
 $x$ such that
$$\lim_{t\to \infty}\frac{a_1+...+a_t}{t}\le \kappa_1 +\varepsilon$$
and
 $?'(x)=0$.}

{\bf Theorem  2.}\,\,\,{\it

1. Let for real irrational $x\in (0,1)$ in continued fraction expansion $
x=[0;a_1,...,a_t,...] $ with $\kappa_2 $ from (\ref{kappa}) the following
inequality be valid:
\begin{equation}
 \liminf_{t\to
\infty}\frac{a_1+...+a_t}{t} >\kappa_2. \label{liminf}
\end{equation}
Then
  $?'(x)$ exists and $?'(x)=0$.

2.For any positive $\varepsilon$ there exists a quadratic irrationality
 $x$ such that
$$\lim_{t\to \infty}\frac{a_1+...+a_t}{t}\ge \kappa_2 -\varepsilon$$
and
 $?'(x)=+\infty$.}

{\bf Theorem 3.}\,\,\,{\it Let in the continued fraction expansion
$x=[0;a_1,...,a_t,...]$ all partial quotients $ a_j$ be bounded by 4. Then
$?'(x) = \infty $.}

We must note that Theorem 3 is not true if we assume that all partial
quotients are bounded by 5.

{\bf Corollary.}\,\,\,{\it The Hausdorff dimension of the set $\{ x:\,\,\,
?'(x) =\infty\}$ is greater than  the Hausdorff dimension of the set ${\cal
F}_4 =\{ x:\,\,\, a_j \le 4 \forall j\}$ which is equal to $0.7889^+$.}

Here the numerical value of Hausdorff dimension for ${\cal F}_4$ is taken from
\cite{Y}. Some resent results on multifractal analysis of the sets associated
with values of $?'(x)$ one can find in the recent paper \cite{ARX}. {\bf
5.\,\,\, The proof of formula  (\ref{MAIN}).}\,\,\, It is sufficient to prove
the unequality
$$ \max_{(r_1,...,r_n)\in
\Omega_{\kappa_2+\eta, n, 1} } \sum_{j=1}^n r_jL_j 
\leq
(L_5-L_4)\eta   .
$$
  By $e_j\in \mathbb{R}^n$ we denote the vector with
all but $j$-th coordinates equal to zero, and with  $j$-th coordinate equal to one.
 The set $\Omega_{\kappa_2+\eta, n, 1}$
is a polytope  lying in the simplex
  $\{r_1,...,r_n:\,\,\, r_j \ge 0, r_1+...+r_n =
1\}$. The vertices of this polytope are points $e_j, 5\le j \le n$ and $e_{i,j}
= \frac{\omega - i}{j-i} e_j +
 \frac{j -\omega}{j-i} e_i,\,\, 1\le i \le 4, 5\le j\le n$.
The linear function $\sum_{j=1}^n r_jL_j$ attend its maximum at a vertex of polytope $\Omega_{\kappa_2+\eta, n, 1}$. Now we must take into
account inequalities
 (\ref{123},\ref{12}).
So we have
$$
\max_{(r_1,...,r_n)\in \Omega_{\kappa_2+\eta, n, t} }  \sum_{j=1}^n r_jL_j = \max\left\{ \max_{1\le i\le 4,\,\, j\ge 5} \left( \left(
\frac{4L_5-5L_4}{L_5 - L_4 } + \frac{jL_i-iL_j}{L_j - L_i }  +\eta \right) \frac{L_j - L_i}{j-i} \right) , L_5\right\}.
$$
But
$$
 \min_{1\le i\le 4,\,\, j\ge 5} \frac{jL_i-iL_j}{L_j - L_i}=
 \frac{5L_4-4L_5}{L_5 - L_4 }
= -\kappa_2
$$
and
$$
 \min_{1\le i\le 4,\,\, j\ge 5,\,\,\,
(i,j) \neq (4,5)
 }
\left( \frac{4L_5-5L_4}{L_5 - L_4 } + \frac{jL_i-iL_j}{L_j - L_i } \right) = \frac{4L_5-5L_4}{L_5 - L_4 }+ \frac{5L_1-L_5}{L_5 - L_1 }
>0.
$$
Hence
$$
\max\left\{ \max_{1\le i\le 4,\,\, j\ge 5} \left( \left( \frac{4L_5-5L_4}{L_5 -
L_4 } + \frac{jL_i-iL_j}{L_j - L_i }  +\eta \right) \frac{L_j - L_i}{j-i}
\right) , L_5\right\} = $$ $$ = \max \left\{ \eta \max_{1\le i\le 4,\,\, j\ge
5} \frac{jL_i-iL_j}{L_j - L_i }, L_5\right\} = \eta (L_5 - L_4).
$$
Formula (\ref{MAIN}) is proved.

 {\bf 6.\,\,\,
One Lemma useful for the proofs of the existence of the derivative of the
Minkowski question mark function.} \,\,\, To prove the existence of the derivative
it is convenient to use the following statement.

{\bf Lemma 2.}\,\,\,{\it For irrational $x$ and  $\delta$ small in absolute
value there exist natural
  $ t = t (x,\delta )$ and $ z\in [1,a_{t+2}+1]$ such that
\begin{equation}
\frac{q_tq_{t-1}}{2^{a_1+...+a_{t+1}+z}} \le \frac{?(x+\delta ) - ?(x)}{\delta }.  \label{DDDD}
\end{equation}
Also    there exist natural
  $ t' = t' (x,\delta )$ and $ z'\in [1,a_{t+2}+1]$ such that
\begin{equation}
\frac{?(x+\delta ) - ?(x)}{\delta } \le \frac{(z'+1)^2q_{t'+1}^2}{2^{a_1+...+a_{t'+1}+z' -4}} \label{DDDD1}
\end{equation}
}

{\bf Proof.}

It is enough to prove Lemma 2 for positive  $\delta$. Define natural
  $n$ such that  $F_n \cap (x,x+\delta ) = \emptyset$,  $F_{n+1} \cap
(x,x+\delta ) = \xi$.
 Let $(x,x+\delta )\subset [\xi^{0},\xi^1]$, where $ \xi^0, \xi^1$ are two
 successive points from the finite set $F_n $. Then $\xi =\xi^0\oplus \xi^1$.
One can easily see that for some natural $t$ will happen $\xi^0 = p_t/q_t$. At
the same time rationals $\xi$ and $\xi^1$ must be among convergent fractions to
$x$ or intermediate fractions to $x$ (intermediate fraction is a fraction of
the form $\frac{p_ta+p_{t-1}}{q_t a+q_{t-1}}, 1\le a< a_{t+1}$).
 In any case, $\xi^1$ has the denominator $\ge q_{t-1}$. Hence
\begin{equation}
\delta \le \frac{1}{q_t q_{t-1} }. \label{w1}
\end{equation}
Define natural $z$ to be minimal such that either  $ \xi_- =
\xi^0\underbrace{\oplus \xi \oplus ... \oplus \xi}_{z} \in (x,\xi) $ or
 $\xi_{+} = \xi^1 \underbrace{\oplus \xi \oplus ... \oplus \xi}_{z} \in (\xi  , x+\delta)$.
 Then
$ \xi_{- -} = \xi^0\underbrace{\oplus \xi \oplus ... \oplus \xi}_{z-1} \le x $
and
 $\xi_{++} = \xi^1 \underbrace{\oplus \xi \oplus ... \oplus \xi}_{z-1} \ge  x+\delta$.
 As points $\xi_{--}<\xi_-<\xi <\xi_+<\xi_{++}$
are successive points from $F_{n+z+1}$
 and $?(x)$ increases, we have
\begin{equation}
\frac{1}{2^{n+z+1}} \le \min \{ \xi_+ -\xi, \xi - \xi_-\}
\le
 ?(x+\delta ) - ?(x) \le ?(\xi_{++}) - ?(\xi_{--}) = \frac{4}{2^{n+z+1}}.
  \label{w2}
\end{equation}

Consider two cases:

(i) \,\,  $\xi_- \in (x,\xi )$.

(ii)  \,\,  $\xi_- \not\in (x,\xi )$ but then $\xi_{+} \in (\xi , x+\delta )$.

 In the case (i) we have $\delta > \xi - \xi_-$. If in addition
(case (i1))
 $ z=1$ then $\xi_- =p/q,
q = z_* q_t + q_{t-1} \le q_{t+1}, 1\le z_* \le a_{t+1},
  \xi =(p - p_t)/(q-q_t)$,
$n+2 = a_1+...+a_{t}+z_*  \le a_1+...+a_{t+1}$
  and
\begin{equation}
\delta > \frac{1}{(q-q_t)q} \ge \frac{1}{(z_*+1)^2q_{t}^2}.
 \label{w3}
\end{equation}
If $ z > 1$ (case (i2))
 then $\xi =p_{t+1}/q_{t+1}$, $ \xi_{--} = p_{t+2}/q_{t+2}, z =
a_{t+2}+1$, $n+1 = a_1+...+a_{t+1}$
 and
\begin{equation}
\delta > \frac{1}{(zq_{t+1}+q_t)q_{t+1}} \ge \frac{1}{(z+1)q_{t+1}^2}.
 \label{w4}
\end{equation}
 In the case (ii) we have
 $ z\le a_{t+2}$,$ \xi =p_{t+1}/q_{t+1}$,
$n+1 = a_1+...+a_{t+1}$. Now we deduce
\begin{equation}
\delta > \xi_+ -\xi \ge \frac{1}{(zq_{t+1}+q^1)q_{t+1}} \ge \frac{1}{(z+1)q_{t+1}^2}
 \label{w5}
\end{equation}
(here $ q^1 < q_{t+1}$ is the denominator of $\xi^1$).

From (\ref{w4},\ref{w5}) and the equalities for $a_1+...+a_{t+1}$   the cases (i2), (ii)
  we get
\begin{equation}
  \delta >   \frac{1}{(z+1)q_{t+1}^2} .
  \label{w55}
\end{equation}

In the cases (i2), (ii) we have $ a_1+...+a_{t+1} -1\le n+1 \le a_1+...+a_{t+1}$. Taking into account (\ref{w1},\ref{w2})and (\ref{w55}) we
 obtain
 $$
\frac{q_tq_{t-1}}{2^{a_1+...+a_{t+1}+z}} \le \frac{?(x+\delta ) - ?(x)}{\delta }  \le \frac{(z+1)q_{t+1}^2}{2^{a_1+...+a_{t+1}+z -4}}
$$
and inequalities (\ref{DDDD},\ref{DDDD1}) follows with $t=t',z=z'$. We should note that the inequality (\ref{DDDD}) also is valid for the case
(i1)as we have $n+2\le a_1+...+a_{t+1}$ and (\ref{w1},\ref{w2}). As for the upper bound in the case (i1) it follows from (\ref{w2},\ref{w3})
with $t' = t-1, $ and $z'=z_*$.

  Lemma 2 is proved.

{\bf 7.\,\,\, The proof of Theorem 1.} \,\,\, The existence of the derivative and its equality to $+\infty$ in the first statement of theorem 1
follows from the lower bound of Lemma 2 as we always have
  $
q_tq_{t-1} \gg \lambda_1^{2t}$ and from the inequality
 $a_1+...+a_{t+1}+a_{t+2}+1 \le \kappa t + o
(t)$ (take into account that $\kappa =
 \limsup_{t\to
\infty}\frac{a_1+...+a_t}{t} <\kappa_1$).

 In order to prove statement 2 of Theorem 1 for small positive
$\eta  >0$ and natural $r$ we define $ q=r^2, m = [r(\kappa_1 -1+\eta)]+1
> r(\kappa_1-1+\eta)$.
Now we  must take the quadratic irrationality
$$
x_r = [0;a_1,...,a_t,...]= [ 0; \overline{ \underbrace{1,..., 1}_{q},  \underbrace{m,..., m}_{r}   } ].
$$
Now we see that
$$
\lim_{t\to \infty}\frac{a_1+...+a_t}{t} = \frac{ q+mr}{q+r} \to \kappa_1 +\eta,\,\,\,\, r\to \infty.
$$
Moreover, taking $ w =\left[\frac{t}{q+r}\right]$ we have
$$
\frac{q_{t+1}(q_{t+1}+q_{t+2})}{2^{a_1+...+a_t}} \le \frac{12 m^3 (k_t(a_1,...,a_t))^2}{2^{a_1+...+a_t}} \le \frac{12
m^32^{2w}\lambda_1^{2wq}\lambda_m^{2wr}}{2^{w(q+rm)}}\le \exp( (-\eta r^2 +O(r\log r))w\log 2).
$$
Here in the exponent the coefficient  before
   $w$
is negative when  $r$ is large enough. Hence the right hand side goes to zero
when   $ t\to \infty$. It means that  $?'(x_r) =0$.

 {\bf 8.\,\,\, The proof of the statement 1 of Theorem 2.} \,\,\,
By Lemma 2 it is sufficient to prove that
   $\frac{q_t^2}{2^{a_1+...+a_t}} \to 0, \,\, t\to
\infty$. Define
 $n$ and $ r_1,...,r_n$
from the condition
 $(a_1,...,a_t) \in W_n (r_1,...,r_n)$. Then  (\ref{PRT}) leads to
$$
\frac{q_t^2}{2^{a_1+...+a_t}} \le \frac{(\mu_n (r_1,...,r_n))^2 }{ 2^{\sum_{j=1}^n jr_j}}
\ll 
\exp \left( 2 \sum_{j=1}^n r_jL_j\right).
$$
From another hand for positive
  $\eta$
small enough  we have the following situation. For all
  $t$ large enough  it is true that
$ n \ge 5$ and $(r_1,...,r_n)\in \Omega_{\kappa_2+\eta,n,t}$. Now we can use
(\ref{MAIN}) and we obtain inequality
$$
\frac{q_t^2}{2^{a_1+...+a_t}} \le \exp \left(  2(L_5-L_4)t\eta\right) \to
0,\,\,\, t\to \infty .
$$
It means that $?'(x) = 0$.

 {\bf 9.\,\,\,The proof of the statement 2 of Theorem 2.}\,\,\,
Take natural numbers $p,q\in \mathbb{N}$ such that $\kappa_2-\varepsilon <
\frac{4p+5q}{p+q}< \kappa_2$. Define
$$
x_{p,q}
 = [0;a_1,...,a_t,...]= [ 0; \overline{ \underbrace{4,..., 4}_{p},  \underbrace{5,..., 5}_{q}   } ].
$$
Obviously,
$$
\lim_{t\to \infty}\frac{a_1+...+a_t}{t} = \frac{4p+5q}{p+q} .
$$
From the other hand
$$
\frac{q_tq_{t-1}}{2^{a_1+...+a_{t+2}}} \ge \left(
 \frac{\lambda_4^{2p}\lambda_5^{2q}}{2^{4p+5q}}\right)^{t+o(t)}=
 \exp (2(pL_4+qL_5)(t+o(t))).
 $$
But $\frac{4p+5q}{p+q}< \kappa_2 = \frac{4L_5-5L_4}{L_5-L_4}$ and hence $
pL_4+qL_5
>0$. So $ \frac{q_tq_{t-1}}{2^{a_1+...+a_{t+2}}} \to \infty$ and $
?'(x_{p,q}) = \infty$.

 {\bf 10.\,\,\, The proof of Theorem  3.} \,\,\,
First of all we see that
\begin{equation}
\min_{ a_i \in \{1,2,3,4\}, a_1+...+a_t = n} \,\,\,\, k_t (a_1,...,a_t) \ge
 \label{MI}
 \end{equation}
$$
   \ge \min\left\{ \min_{ a_i \in \{1,4\},
 a_1+...+a_t = n-3} \,\, k_t (a_1,...,a_t),
 \min_{ a_i \in \{1,4\}, a_1+...+a_t = n-2} \,\, k_t (a_1,...,a_t),
\min_{ a_i \in \{1,4\}, a_1+...+a_t = n} \,\, k_t (a_1,...,a_t)
 \right\}.
  $$
 In order to do this we
note that for two elements $a,b$ with other elements fixed
$$
k_t (...,a,...,b,...) = Mab+Na+Kb +P .
$$
Here positive $M,N,K,P$ do not depend on $a,b$. Then if the sum $a+b =\tau$ is fixed we have
$$
k_t (...,a,...,b,...) = Ma(\tau - a)+Na+K(\tau - a) +P = -Ma^2+ (M\tau+N - K) a -K\tau +P.
$$
So for $ a,b >1$ we can say that
$$
 k_t (...,a,...,b,...)\ge \min \{ k_t (...,a- 1,...,b+1,...), k_t (...,a+1,...,b-1,...)\}.
$$
Hence, we can replace a pair  $2,3$ of partial quotients by $1,4$ and the continuant becomes smaller. Also we can replace any pair  $2,2$ of
partial quotients by $1,3$ and the continuant becomes smaller. Also we can replace any pair  $3,3$ of partial quotients by $2,4$ and the
continuant becomes smaller. This procedure enables one to replace the set $\{ (a_1,...,a_t) :\,\,\, a_i \in \{1,2,3,4\}, a_1+...+a_t = n\}$ in
the left hand side of (\ref{MI}) by the set $\{ (a_1,...,a_t) :\,\,\, a_i \in \{1,2,3,4\}, a_1+...+a_t = n,\,\,\, \#\{ a_i = 3 \}+\#\{ a_i = 2 \}
\le 1 \}$. Now the inequality (\ref{MI}) follows.

From another hand as all partial quotients are bounded by $4$ we have
$$
k_{t_1+t_2} (a_1,...,a_{t_1},a_1,...,a_{t_2}) 
\ge 
\left(1+\varepsilon \right) k_{t_1} (a_1,...,a_{t_1}) k_{t_2} (a_1,...,a_{t_2}), 
$$
where $\varepsilon$ is some relatively small positive real constant.
Now from the last formulas and  (\ref{MI}) it follows that it is sufficient to prove  that for every large $n$ the following inequality is valid
  \begin{equation}
\min_{a_1+...+a_t=n,a_j \in \{ 1, 4\}} k_t (
   a_1,...,a_t
   )
 \ge (\sqrt{2})^n
 \label{sqrt}
 \end{equation}
 (here minimum is taken over all $t$-tuples
$a_1,...,a_t$ such that
 $a_1+...+a_t=n$ and $ a_j \in \{ 1, 4\} $).
 This can be easy verified by induction in $n$.
The base of induction for $n = 23, 24$ is checked by computer by MAPLE (the program is given in section 10). By the Sylvester theorem any
natural number $t$ greater than $505 = 23\times 24 - 23 -24$ can be expressed in the form $ t = 23x+24y$ with nonnegative  integers $x,y$. Hence
for $t \ge 506$ we have
$$
k_t (
   a_1,...,a_t
   )
   \ge
   \prod_{1\le j\le x}
   k_{23} (
   a_1^{(j)},...,a_{23}^{(j)}
   )
   \prod_{1\le j\le y}
   k_{24} (
   b_1^{(j)},...,b_{24}^{(j)}
   )
   $$
(here $ (
   a_1,...,a_t
   ) =
(
   a_1^{(1)},...,a_{23}^{(1)},...,
      a_1^{(x)},...,a_{23}^{(x)},
      b_1^{(1)},...,b_{24}^{(1)},...,
      b_1^{(y)},...,b_{24}^{(y)}
      )$).

   Now (\ref{sqrt}) follows from the base of induction for $n = 23, 24$.
   Theorem 3 is proved.

 {\bf 11.\,\,\,
MAPLE program for verifying the inequalities for $n=23,24$.}\,\,\, Here is the program for $ n=23$. The program for $ n = 24$ is quite similar.

\vskip+1.0cm

 for
 $ a_1$ from 1 by 3 to 4 do

 for $a_2$ from 1 by 3 to 4 do

 for $a_3$ from 1 by 3 to 4 do

 for $a_4$ from 1 by 3 to 4 do

 for $a_5$ from 1 by 3 to 4 do

 for $a_6$ from 1 by 3 to 4 do

 for $a_7$ from 1 by 3 to 4 do

 for $a_8$ from 1 by 3 to 4 do

 for $a_9$ from 1 by 3 to 4 do

 for $a_{10}$ from 1 by 3 to 4 do

 for $a_{11}$ from 1 by 3 to 4 do

 for $a_{12}$ from 1 by 3 to 4 do

 for $a_{13}$ from 1 by 3 to 4 do

 for $a_{14}$ from 1 by 3 to 4 do

 for $a_{15}$ from 1 by 3 to 4 do

 for $a_{16}$ from 1 by 3 to 4 do

 for $a_{17}$ from 1 by 3 to 4 do

 for $a_{18}$ from 1 by 3 to 4 do

 for $a_{19}$ from 1 by 3 to 4 do

 for $a_{20}$ from 1 by 3 to 4 do

 for $a_{21}$ from 1 by 3 to 4 do

 for $a_{22}$ from 1 by 3 to 4 do

 for $a_{23}$ from 1 by 3 to 4 do

 $k_1:=a_1$;

 $k_2:=a_2*k_1+1$;

 $k_3:=a_3*k_2+k_1$;

 $k_4:=a_4*k_3+k_2$;

 $k_5:=a_5*k_4+k_3$;

 $k_6:=a_6*k_5+k_4$;

 $k_7:=a_7*k_6+k_5$;

 $k_8:=a_8*k_7+k_6$;

 $k_9:=a_9*k_8+k_7$;

 $k_{10}:=a_{10}*k_9+k_8$;

 $k_{11}:=a_{11}*k_{10}+k_9$;

 $k_{12}:=a_{12}*k_{11}+k_{10}$;

 $k_{13}:=a_{13}*k_{12}+k_{11}$;

 $k_{14}:=a_{14}*k_{13}+k_{12}$;

 $k_{15}:=a_{15}*k_{14}+k_{13}$;

 $k_{16}:=a_{16}*k_{15}+k_{14}$;

 $k_{17}:=a_{17}*k_{16}+k_{15}$;

 $k_{18}:=a_{18}*k_{17}+k_{16}$;

 $k_{19}:=a_{19}*k_{18}+k_{17}$;

 $k_{20}:=a_{20}*k_{19}+k_{18}$;

 $k_{21}:=a_{21}*k_{20}+k_{19}$;

 $k_{22}:=a_{22}*k_{21}+k_{20}$;

 $k_{23}:=a_{23}*k_{22}+k_{21}$;

 $e_{23}:=
2^{(a_1+a_2+a_3+a_4+a_5+a_6+a_7+a_8+a_9+a_{10}+a_{11}+a_{12}+a_{13}+a_{14}+a_{15}+a_{16}+a_{17}+ a_{18}+a_{19}+a_{20} +a_{21}+a_{22}+a_{23})}$;

 if($(k_{23})^2<e_{23}$) then

 print($a_1,a_2,a_3,a_4,a_5,a_6,a_7,a_8,a_9,a_{10}, a_{11},a_{12},a_{13},a_{14},a_{15},a_{16},a_{17},a_{18},a_{19},a_{20},a_{21},a_{22},a_{23}$);

 end if;

  end do;

 end do;

 end do;

 end do;

 end do;

 end do;

 end do;

 end do;

 end do;

 end do;

 end do;

 end do;

 end do;

 end do;

 end do;

 end do;

 end do;

 end do;

 end do;

 end do;

 end do;

 end do;

 end do;

\vskip+1.0cm

Authors:

\vskip+1.0cm

Moshchevitin Nikolai G.,
e-mail: moshchevitin@rambler.ru
\vskip+1.0cm

Dushistova Anna A.,
e-mail: anchatnik@bk.ru

\end{document}